\documentclass[12pt]{amsart}
\usepackage{amsmath, amstext, amsbsy, amssymb}

\newtheorem{corollary}{Corollary}[section]

\newtheorem{theorem}{Theorem}[section]
\newtheorem{proposition}{Proposition}[section]

\numberwithin{equation}{section}

\newcommand{\ga}{\gamma}
\newcommand{\la}{\lambda}

\newcommand{\ula}{\underline{\lambda}}
\newcommand{\umu}{\underline{\mu}}
\newcommand{\unu}{\underline{\nu}}
\newcommand{\urh}{\underline{\rho}}

\newcommand{\uth}{\underline{\theta}}
\newcommand{\uT}{\underline{T}}
\newcommand{\uw}{\underline{w}}
\newcommand{\LQ}{\Lambda_{\mathbb Q}}

\newcommand{\LZ}{\Lambda_{\mathbb Z}}
\newcommand{\LQr}{\LQ^{\otimes r}}
\newcommand{\LZr}{\LZ^{\otimes r}}

\begin{document}

\title[Wreath Product Symmetric Functions]
{Wreath product symmetric functions}
\author{Frank Ingram}
\address{FI: Mathematics Department, Winston-Salem State University,
Winston-Salem, NC 27110, USA}
\email{ingramfr@wssu.edu}
\author{Naihuan Jing$^*$}
\address{NJ: Department of Mathematics,
   North Carolina State University,
   Raleigh, NC 27695-8205, USA}
\email{jing@math.ncsu.edu}
\author{Ernie Stitzinger}
\address{ES: Department of Mathematics,
   North Carolina State University,
   Raleigh, NC 27695-8205, USA}
\email{stitz@math.ncsu.edu}
\thanks{Jing greatly acknowledge the support of NSA grant
H98230-06-1-0083 and  Alexander von Humboldt Fellowship.
\thanks{$^*$ \it Corresponding author}}

\begin{abstract}
We systematically study wreath product Schur functions and give a
combinatorial construction using colored partitions and tableaux.
The Pieri rule and the Littlewood-Richardson rule are studied. We
also discuss the connection with representations of generalized
symmetric groups.
\end{abstract}
\maketitle


\section{Introduction}

Since Jacobi and Frobenius, symmetric functions have played a
fundamental role in representations theory and number theory. The
character formula of the symmetric group in terms of Schur functions
and power sums is undoubtedly one of the most beautiful chapters in
group theory and algebraic combinatorics. Schur developed Frobenius
theory further using his rational homogeneous functions, and again
symmetric functions were the most important techniques in his work
on generalizations of Frobenius character theory.

In Macdonald's book \cite{M}, one can see other interesting examples
of applications of symmetric functions into group theory. Among
them, Specht's generalization \cite{Sp} of Schur functions deserves
special mention. Let $G$ be any finite group, the semidirect product
$G^n\ltimes S_n$ denoted by $G\sim S_n$ is called the wreath product
of the finite group $G$ and the symmetric group $S_n$. The special
case $1\sim S_n$ is exactly the symmetric group $S_n$, while
$\mathbb Z_2\sim S_n$ is the hyperoctahedral group. Recently the
wreath products
are used in the generalized McKay
correspondence \cite{FJW} to realize the affine Lie algebras of
simply laced types. The wreath products have also
been studied in various interesting contexts \cite{H, HH, MR}.
However the basic formulation of the
characters remains complicated and formidable.

It is well-known that
Schur functions can also be defined in terms of Young tableaux in
Stanley's work \cite{St}. One purpose of this paper is to
give a simpler and more elementary
description of wreath product symmetric functions
using colored tableaux. We achieve this by first considering tensor
product of Schur functions and then pass the results to the case of
wreath product using linear transformations. In our simplified
approach one can see the structures of wreath product
symmetric functions easily.
In the same spirit we also obtain the description of
the Littlewood-Richardson rule for the wreath product Schur functions.
We deviate slightly from Macdonald's inner product by requiring that
the power sum associated with irreducible characters of $G$ to be
invariant under the complex conjugation, which will facilitate us to
pass the results from tensor product to wreath product easily.

The paper is organized into five parts. After introduction (section
one) we give a preliminary review of the main results concerning
Schur functions to set the notations and for later usage in Section
two. In Section three we studied tensor products of symmetric
functions and formulate several bases modeled on one variable case.
The notion of colored partitions and tableaux is the fundamental
combinatorial objects for us to formulate results in tensor
products. In Section four wreath product Schur functions are
introduced by two sets of power sum symmetric functions indexed by
colored partitions, which correspond to conjugacy classes and
irreducible characters of $G\sim S_n$. We show how to pass from
tensor products to wreath products by carefully treating the inner
product. In section five we give generalized Littlewood-Richardson
rule and discuss their special cases, and finally we discuss the
relationship between symmetric functions and characters of wreath
products of symmetric groups.


\section{Symmetric functions}\label{S:symm}

In this section we review some basic materials on symmetric
functions that will be used in our later discussion. Almost all
statements are standard and can be found in \cite{M} or in a
slightly different notation in \cite{St}.

A \textit{partition} $\la$ of an non-negative integer $n$, denoted
as $\la\vdash n$, is an integral decomposition of $n$ arranged in
decreasing order: $n=\la_1+\ldots+\la_l$, $\la_1\geq\ldots\geq
\la_l\geq 0$. If the last part $\la_l\neq 0$, the partition has
\textit{length} $l(\la)=l$, and we also denote by $|\la|=n$. One
sometimes uses the notation $\la=(1^{m_1}2^{m_2}\cdots)$ if there
are $m_i$ parts equal to $i$ in $\la$.

The \textit{Ferrer diagram} associated with the partition $\la\vdash
n$ with $l(\la)=l$ consists of $l$ layers of $n$ boxes aligned to
the left from the top to the bottom. The partition corresponding to
the flip of the Ferrer diagram of $\la$ along the
northwest-to-southeast axis is called the \textit{conjugate
partition} $\la'$. For example the Ferrer diagram $(43221)$ and its
conjugate diagram are

\begin{table}[!htpb]
\begin{center}
\begin{tabular}{lll}
\framebox[.28in]{\rule{0in}{.2in}}\framebox[.28in]{\rule{0in}{.2in}}\framebox[.28in]{\rule{0in}{.2in}}\framebox[.28in]{\rule{0in}{.2in}}
&&\framebox[.28in]{\rule{0in}{.2in}}\framebox[.28in]{\rule{0in}{.2in}}\framebox[.28in]{\rule{0in}{.2in}}\framebox[.28in]{\rule{0in}{.2in}}\framebox[.28in]{\rule{0in}{.2in}}\\
\framebox[.28in]{\rule{0in}{.2in}}\framebox[.28in]{\rule{0in}{.2in}}\framebox[.28in]{\rule{0in}{.2in}}
&&\framebox[.28in]{\rule{0in}{.2in}}\framebox[.28in]{\rule{0in}{.2in}}\framebox[.28in]{\rule{0in}{.2in}}\framebox[.28in]{\rule{0in}{.2in}}\\
\framebox[.28in]{\rule{0in}{.2in}}\framebox[.28in]{\rule{0in}{.2in}}&&\framebox[.28in]{\rule{0in}{.2in}}\framebox[.28in]{\rule{0in}{.2in}}\\
\framebox[.28in]{\rule{0in}{.2in}}\framebox[.28in]{\rule{0in}{.2in}}&&\framebox[.28in]{\rule{0in}{.2in}}\\
\framebox[.28in]{\rule{0in}{.2in}}
\end{tabular}
\caption{$\lambda=(4,3,2,2,1)$  and $\lambda^{\prime}=(5,4,2,1)$}
\end{center}
\end{table}
There are two partial orders in the set of partitions. The first one
is the \textit{dominance order} $\geq$. Let $\la$ and $\mu$ be two
partitions, $\la\geq \mu$ means that $\la_1\geq\mu_1$,
$\la_1+\la_2\geq\mu_1+\mu_2$, etc. The second one is the
\textit{containment order} $\supset$. If all parts of $\la$ are
larger than corresponding parts of $\mu$, we say that $\la$ contains
$\mu$ and denote by $\la\supset\mu$.

Let $\LZ$ be the ring of symmetric functions in the $x_i (i\in
\mathbb N)$ over $\mathbb Z$. For each non-negative integer $n$, the
\textit{elementary symmetric function}  $e_n$ is defined by
$$e_n=\sum_{i_1<\cdots<i_n}x_{i_1}\cdots x_{i_n},$$
where the sum is taken formally, or it is in the sense of the
inverse direct limit determined by the relations: $e_m(x_1, \ldots,
x_n, 0, \ldots, 0)=e_n(x_1, \ldots, x_n)$ whenever $m\geq n$.

For each partition $\la$ with $l(\la)\leq n$ the {\it monomial
symmetric function} $m_{\la}$ is defined by $m_{\la}(x_1,\ldots,
x_n)=x_1^{\la_1}\cdots x_{\la_n}^{\la_n}+\mbox{distinct
permutations}$, where we take $\la_{l+1}=\cdots=\la_n=0$ if
necessary.

For each partition $\la$ we define
\begin{align*}
e_{\la}&=e_{\la_1}\cdots e_{\la_l},\\
m_{\la}&=m_{\la_1}\cdots m_{\la_l}.
\end{align*}

The \textit{homogeneous symmetric function} $h_n$ is defined by
$h_n=\sum_{|\la|=n}m_{\la}$. For any partition $\la$ we define
$h_{\la}=h_{\la_1}\cdots h_{\la_l}$

As a graded vector space we have
\begin{equation}
\LZ=\bigoplus_{\la}\mathbb Ze_{\la} =\bigoplus_{\la}\mathbb Zm_{\la}
=\bigoplus_{\la}\mathbb Zh_{\la},
\end{equation}
where the sums are over all partitions. Moreover $\LZ=\mathbb Z[e_1,
e_2, \ldots]=\mathbb Z[h_1, h_2, \ldots]$ as polynomials rings, and
$\{e_n\}$ and $\{h_n\}$ are algebraically independent over $\mathbb
Z$. Subsequently the sets $\{e_{\la}\}$, $\{m_{\la}\}$ and
$\{h_{\la}\}$ are all $\mathbb Z$-linear bases of $\LZ$.

For each non-negative integer $n$, the \textit{power sum} symmetric
function $p_n=\sum_ix_i^n$ is an element in $\LQ$, and we define for
each partition $\la$
\begin{equation}
p_{\la}=p_{\la_1}\cdots p_{\la_l}.
\end{equation}
It is well-known that the $\{p_{\la}\}$ forms a $\mathbb Q$-linear
basis in $\LQ$, and subsequently $\LQ=\mathbb Q[p_1, p_2, \ldots]$.
Moreover the $p_n$ are algebraically independent over $\mathbb Q$.

The space $\LZ$ has a $\mathbb Z$-valued bilinear form $<\, ,\, >$
defined by
\begin{equation}
<h_{\la}, m_{\mu}>=\delta_{\la\mu}.
\end{equation}
Equivalently the symmetric bilinear form can also be defined by
\begin{equation}
<p_{\la}, p_{\mu}>=\delta_{\la\mu}z_{\mu},
\end{equation}
where $z_{\la}=\prod_ii^{m_i}m_i!$. The distinguished orthonormal
basis--\textit{Schur functions} are given by the following
triangular relations \cite{M}:
\begin{align*}
&s_{\la}=m_{\la}+\sum_{\mu<\la}c_{\la\mu}m_{\mu}, \qquad \mbox{for
some $c_{\la\mu}\in \mathbb Z$} \\
&<s_{\la}, s_{\mu}>=\delta_{\la\mu}
\end{align*}

Let $\lambda$ and $\mu$ be partitions such that $\mu\subset\lambda$.
The \textit{skew diagram} $\lambda/\mu$ is the set-theoretic
difference $\lambda-\mu$. A maximal connected component of $\la/\mu$
is called a connected component. Each connected component itself is
a skew diagram. For example the skew Ferrer diagram $\lambda /\mu$
=(6,4,2,2)/(4,2) has three components and is shown by
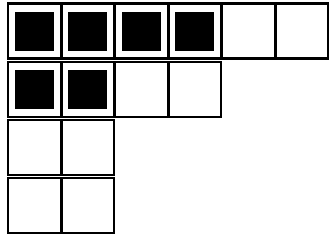
\begin{table}[!hpb]
\begin{tabular}{llllll}
\framebox[.28in]{\rule{.2in}{.2in}}\framebox[.28in]{\rule{0.2in}{.2in}}\framebox[.28in]{\rule{0.2in}{.2in}}\framebox[.28in]{\rule{0.2in}{.2in}}\framebox[.28in]{\rule{0in}{.2in}}\framebox[.28in]{\rule{0in}{.2in}}\\
\framebox[.28in]{\rule{0.2in}{.2in}}\framebox[.28in]{\rule{0.2in}{.2in}}\framebox[.28in]{\rule{0in}{.2in}}\framebox[.28in]{\rule{0in}{.2in}}\\
\framebox[.28in]{\rule{0in}{.2in}}\framebox[.28in]{\rule{0in}{.2in}}\\
\framebox[.28in]{\rule{0in}{.2in}}\framebox[.28in]{\rule{0in}{.2in}}\\

\end{tabular}
\end{table}

A \textit{horizontal strip} is a  skew Ferrer diagram $\lambda /\mu
$ with no two boxes in the same column, i.e., $\la'_i-\mu_i'\leq 1$.
A \textit{vertical strip} is skew Ferrer diagram $\lambda /\mu $
with no two boxes in the same row, i.e., i.e., $\la_i-\mu_i\leq 1$


The multiplication of two Schur functions is again a symmetric
function and thus can be expressed as a linear combination of Schur
functions:
\begin{equation}
s_{\la}s_{\mu}=\sum_{\nu}c_{\la\mu}^{\nu}s_{\nu}, \qquad
c_{\la\mu}^{\nu}\in\mathbb Z_+.
\end{equation}
The structure constants $c_{\la\mu}^{\nu}$ are called the
\textit{Littlewood-Richardson coefficients}. It is known that
$c_{\la\mu}^{\nu}\neq 0$ unless $|\nu|=|\la|+|\mu|$ and $\nu\supset
\la, \mu$.

The \textit{skew Schur function} $s_{\la/\mu}\in\LZ$ is then defined
by
\begin{equation}
s_{\la/\mu}=\sum_{\nu}c_{\mu\nu}^{\la}s_{\nu}.
\end{equation}
It is clear that $deg(s_{\la/\mu})=|\la|-|\mu|$ and $s_{\la/\mu}=0$
unless $\mu\subset \la$, i.e. $s_{\la/\mu}$ is defined for the skew
Ferrer diagram $\la/\mu$.

Let $x=(x_1, x_2, \ldots)$ and $y=(y_1, y_2, \ldots)$ be two sets of
variables.
\begin{proposition} \cite{M} \label{schur}
(a) If the connected components of the skew diagram $\la/\mu$ are
$\theta_i$, then $s_{\la/\mu}(x)=\prod_is_{\theta_i}(x)$. In
particular, $s_{\la/\mu}(x)=0$ if $\mu \nsubseteq \la$.

 (b) The skew Schur symmetric function
$s_{\la}(x, y)$ satisfies
\begin{equation}
s_{\la/\mu}(x, y)=\sum_{\nu} s_{\la/\nu}(x)s_{\nu/\mu}(y),
\end{equation}
where the sum runs through all the partition $\nu$ such that
$\la\supset \nu \supset \mu$.

(c) In general, the skew Schur function $s_{\la/\mu}(x^{(1)},
\ldots, x^{(n)})$ can be written as
\begin{equation}
s_{\la/\mu}(x^{(1)}, \ldots, x^{(n)})=\sum_{(\nu)}\prod_{i=1}^n
s_{\nu^{(i)}/\nu^{(i-1)}}(x^{(i)}),
\end{equation}
where the sum runs through all sequences of partitions
$(\nu)=(\nu^{(n)}, \ldots, \nu^{(0)})$ such that
$\mu=\nu^{(0)}\subset \nu^{(1)} \subset \ldots \subset\nu^{(n)}
=\la$.
\end{proposition}

If there is only one variable, the Schur function can be easily
computed.

\begin{proposition} \cite{M} \label{skew} (a) We have $s_{\la}(x)=0$ when $l(\la)>1$.
When $\la/\mu$ is a horizontal strip,
$s_{\la/\mu}(x)=x^{|\la|-|\mu|}$.

(b) More generally $s_{\la/\mu}(x_1, \ldots, x_n)=0$ unless
$\la_i'-\mu_i'\leq n$ for each $i$.
\end{proposition}

In order to describe Schur functions combinatorially we introduce
the concept of Young tableaux.

 For a diagram $\la$, a \textit{semistandard tableau} $T$ of shape $\la$ is an insertion
 of natural numbers $1, 2, \ldots $ into $\la$ such that the rows weakly increase
 and the columns strictly increase.
The \textit{content} $\mu$ of a tableau T is the composition
$\mu=\left(  \mu_{1},\mu_{2},\ldots,\mu_{n}\right)$, where $\mu_{i}$
equals the number of $i$'s in T. We will write $x^T=x_1^{\mu_1}
\cdots x_n^{\mu_n}$. In many situations we want to restrict the
largest possible integer $n$ inserted into the diagram $\la$, and
the notation $x^T$ makes it clear when one uses $x=(x_1, \ldots,
x_n)$.

A \textit{semistandard Young Tableau}, denoted by $T_{\lambda\mu}$,
is a tableau T having shape $\lambda$ and content $\mu$. For
example, see Table 2. 

\begin{table}[!hpb]
\begin{center}
\begin{tabular}{l}
\framebox[.28in]{\rule{0in}{.2in}{1}}\framebox[.28in]{\rule{0in}{.2in}{1}}\framebox[.28in]{\rule{0in}{.2in}{1}}\framebox[.28in]{\rule{0in}{.2in}{4}}\\\nolinebreak
\framebox[.28in]{\rule{0in}{.2in}{2}}\framebox[.28in]{\rule{0in}{.2in}{3}}\framebox[.28in]{\rule{0in}{.2in}{3}}\\
\framebox[.28in]{\rule{0in}{.2in}{3}}\framebox[.28in]{\rule{0in}{.2in}{4}}\\
\framebox[.28in]{\rule{0in}{.2in}{5}}\framebox[.28in]{\rule{0in}{.2in}{5}}\\
\framebox[.28in]{\rule{0in}{.2in}{6}}
\end{tabular}
\caption{$T_{\lambda\mu}$ with $\lambda=(4,3,2,2,1)$ and
$\mu=(3,1,3,2,2,1)$}
\end{center}
\end{table}

The \textit{Kostka number} $K_{\lambda\mu}$ is the number of
semistandard tableaux of shape $\lambda$ and content $\mu$.

Let $\lambda $ and $\mu $ be partitions  such that $\lambda
\supseteq $ $\mu $ (i.e. $\lambda _{i}\geq \mu _{i}$ for all $i$). A
\textit{\ semistandard Young tableau of skew shape $\lambda /\mu $}
with content $\nu$, is a diagram of shape $\lambda /\mu $ \ whose
boxes have been filled with $|\nu|$ positive integers that are
weakly increasing in every row and strictly increasing in every
column. For example, the skew tableau of shape $\lambda /\mu$
=(6,5,4,3)/(4,2) with content $\nu=(2,3,3,0,3,1)$ is shown here

\begin{table}[!hpb]
\begin{tabular}{llllll}
\framebox[.28in]{\rule{0.2in}{.2in}}\framebox[.28in]{\rule{0.2in}{.2in}}\framebox[.28in]{\rule{0.2in}{.2in}}\framebox[.28in]{\rule{0.2in}{.2in}}\framebox[.28in]{\rule{0in}{.2in}{2}}\framebox[.28in]{\rule{0in}{.2in}{2}}\\
\framebox[.28in]{\rule{0.2in}{.2in}}\framebox[.28in]{\rule{0.2in}{.2in}}\framebox[.28in]{\rule{0in}{.2in}{1}}\framebox[.28in]{\rule{0in}{.2in}{1}}\framebox[.28in]{\rule{0in}{.2in}{5}}\\
\framebox[.28in]{\rule{0in}{.2in}{2}}\framebox[.28in]{\rule{0in}{.2in}{3}}\framebox[.28in]{\rule{0in}{.2in}{3}}\framebox[.28in]{\rule{0in}{.2in}{3}}\\
\framebox[.28in]{\rule{0in}{.2in}{5}}\framebox[.28in]{\rule{0in}{.2in}{5}}\framebox[.28in]{\rule{0in}{.2in}{6}}\\
\end{tabular}
\caption{tableau $(6,5,4,3)/(4,2)$ with content $(2,3,3,0,3,1)$}
\end{table}
which we will also write as the array

\begin{center}
$\begin{array}{cccccc}
&  &  &  & 2 & 2 \\
&  & 1 & 1 & 5 & \\
2 & 3 & 3 & 3 &  &  \\
4& 5 & 6 &  &  &
\end{array}$
\end{center}
Similarly the \textit{Kostka number} $K_{\lambda/\mu, \nu}$ is the
number of semistandard tableaux of shape $\lambda/\mu$ and content
$\nu$.

Combining Propositions \ref{skew} and \ref{schur} one easily gets
the following result

\begin{proposition} The Schur function $s_{\la/\mu}(x)$ can be
expressed as a summation of all monomials $x^T$ attached to Young
tableaux of shape $\la/\mu$:
\begin{equation}
s_{\la/\mu}(x)=\sum_T x^T,
\end{equation}
where the sum runs through all skew semistandard Young tableaux of
shape $\la/\mu$.
\end{proposition}
\begin{proof} In Proposition \ref{skew}(c) we take $x^{(i)}=x_i$ and using
Proposition \ref{schur} it follows that
\begin{align*}
s_{\la/\mu}(x_1, \ldots, x_n)&=\sum_{(\nu)}\prod_{i=1}^n
s_{\nu^{(i)}/\nu^{(i-1)}}(x_i)\\
&=\sum_{(\nu)}\prod_{i=1}^n x_i^{|\nu^{(i)}|-|\nu^{(i-1)}|},
\end{align*}
where the sum runs through all sequences of partitions
$(\nu)=(\nu^{(n)}, \ldots, \nu^{(0)})$ such that
$\mu=\nu^{(0)}\subset \nu^{(1)} \subset \ldots \subset\nu^{(n)}
=\la$ and $\nu^{(i)}/\nu^{(i-1)}$ are horizontal strips. For each
horizontal strip $\nu^{(i)}/\nu^{(i-1)}$, we assign $i$ in the
diagram $\la/\mu$ where the horizontal strip occupies and thus
obtain a semistandard Young tableau of shape $\la/\mu$ with content
$T=(|\nu^{(1)}|-|\nu^{(0)}|, |\nu^{(2)}|-|\nu^{(1)}|, \ldots,
|\nu^{(n)}|-|\nu^{(n-1)}|)$. In this way the sequence of such
partitions gives rise a semistandard Young tableau of shape
$\la/\mu$ with content $T$. Conversely given a semistandard Young
tableau of shape $\la/\mu$ with content $T$ $(|T|=|\la|-|\mu|)$, let
$\nu^{(i)}$ be the subdiagram consisting of the diagram $\mu$
together with entries numbered $1, \ldots, i$. Then the skew diagram
$\nu^{(i)}/\nu^{(i-1)}$ is horizontal and gives rise a sequence of
partitions $(\nu)=(\nu^{(n)}, \ldots, \nu^{(0)})$ such that
$\mu=\nu^{(0)}\subset \nu^{(1)} \subset \ldots \subset\nu^{(n)}
=\la$ and $\nu^{(i)}/\nu^{(i-1)}$ are horizontal strips. Therefore
we have
$$
s_{\la/\mu}(x)=\sum_{T}x^T,
$$
where $T$ runs through Young tableaux of shape $\la/\mu$.
\end{proof}

\begin{corollary} We have $\displaystyle s_{\la/\mu}=\sum_{\nu}K_{\la/\mu,
\nu}m_{\nu}$.
\end{corollary}

\section{Tensor product of symmetric functions} \label{S:tensor}

Fix a natural number $r$ we consider the tensor product
$\LZ^{\otimes r}$. Clearly the tensor products of basis elements
from a given basis of $\LZ$ will form a basis for the tensor product
space $\LZr$.

To parametrize the basis elements, we introduce the notion of
\textit{colored partitions}. A colored partition $\ula$ is a
partition-valued function: $\ula=(\la^{(0)}, \ldots, \la^{(r-1)})$,
where $\la^{(i)}$ are partitions. Here we intuitively color the
$i$th diagram $\la^{(i)}$ by the color $i$. Also we set $I=\{0, 1,
\ldots, r-1\}$, the set of the indices or \textit{colors}. If we
want to specify the number of colors or partitions inside $\ula$, we
also say that $\ula$ is a $r$-colored partition. We denote by
$|\ula|$ the sum of all weights: $|\ula|=|\la^{(0)}|+\cdots
+|\la^{(r-1)}|$. Thus we will also use the terms such as the
\textit{colored Ferrer diagrams, colored skew Ferrer diagrams} etc.

The dominance order $\leq$ and containment order $\subset$ can be
extended to colored partitions as follows. For two colored
partitions $\ula$ and $\umu$, $\ula\leq\umu$ means that
$\la^{(i)}\leq \mu^{(i)}$ for $i\in I$. Similarly $\ula\subset\umu$
means that $\la^{(i)}\subset\mu^{(i)}$ for $i\in I$.

Let $u_{\la}$ be any element from our bases $\{e_{\la}\}$,
$\{m_{\la}\}$, $\{h_{\la}\}$, $\{p_{\la}\}$ or $\{s_{\la}\}$ of
$\LZ$ in Section \ref{S:symm}. Let $x^{(i)}=(x_{i,1}, x_{i, 2},
\ldots)$ be the variables in the $i$th ring $\LZ$ ($i\in I$), and
for any $r$-colored partition $\ula$ we define
\begin{equation}
u_{\ula}=u_{\la^{(0)}}(x^{(0)})\cdots u_{\la^{(r-1)}}(x^{(r-1)}).
\end{equation}
Then the $u_{\ula}$ forms a basis for the tensor product space
$\LQr$. We also define similarly
$$
z_{\ula}=z_{\la^{(0)}}\cdots z_{\la^{(r-1)}}.
$$

We canonically extend the bilinear scalar product of $\LZ$ to $\LZr$
Namely we define on $\LZr$ by
$$<u_1\otimes\cdots\otimes u_r, v_1\otimes\cdots\otimes v_r>=<u_1,
v_1>\cdots <u_1, v_1>
$$
for any two sets of elements $\{u_i\}$, $\{v_i\} \subset \LZ$ and
extend bilinearly. Then we have the following result.

\begin{theorem} (a) The tensor product ring $\LZr$ of symmetric functions
has the following four sets of linear bases: $\{e_{\ula}\}$,
$\{h_{\ula}\}$, $\{m_{\ula}\}$ and $\{s_{\ula}\}$. The set
$\{p_{\ula}\}$ is a $\mathbb Q$-linear basis for $\LQr$. Namely we
have
\begin{align*}
\LZr&=\bigoplus_{\ula}\mathbb Ze_{\ula} =\bigoplus_{\ula}\mathbb
Zm_{\ula} =\bigoplus_{\ula}\mathbb Zh_{\ula}=\bigoplus_{\ula}\mathbb Zs_{\ula},\\
\LQr&=\bigoplus_{\ula}\mathbb Qp_{\ula},
\end{align*}
where the sums are over all colored partitions. Moreover
\begin{align*}
\LZ&=\bigotimes_{i=0}^{r-1} \mathbb Z[e_{i,1}, e_{i,2}, \ldots]=\bigotimes_{i=0}^{r-1}\mathbb Z[h_{i,1}, h_{i,2}, \ldots]\\
&\simeq\mathbb Z[e_{1}, e_{2}, \ldots]^{\otimes r}\simeq\mathbb
Z[h_{1}, h_{2}, \ldots]^{\otimes r}
\end{align*}
as polynomials rings, and $\{e_n\}$ and $\{h_n\}$ are
algebraically independent over $\mathbb Z$. Subsequently the sets
$\{e_{\la}\}$, $\{m_{\la}\}$ and $\{h_{\la}\}$ are all $\mathbb
Z$-linear bases of $\LZ$.

(b) The bases $\{h_{\ula}\}$ and $\{m_{\ula}\}$ are dual under the
scalar product:
$$
<h_{\ula}, m_{\ula}>=\delta_{\ula\umu}.
$$

 (c) The power sum symmetric functions are orthogonal:
$$
<p_{\ula}, p_{\umu}>=\delta_{\ula\umu}z_{\ula}.
$$

(d) The Schur functions $s_{\ula}$ are orthonormal and uniquely
determined by the triangular relations:
\begin{align*}
&s_{\ula}=m_{\ula}+\sum_{\mu<\la}c_{\ula\umu}m_{\umu}, \qquad
\mbox{for
some $c_{\ula\umu}\in \mathbb Z$} \\
&<s_{\ula}, s_{\umu}>=\delta_{\ula\umu}.
\end{align*}
\end{theorem}
\begin{proof} All statements are trivial when one uses the tensor product structure
$\LZr$ or $\LQr$. For example, one has
\begin{align*}
s_{\ula}&=s_{\la^{(0)}}\cdots s_{\la^{(r-1)}}\\
&=(\sum_{\mu^{(0)}\leq
\la^{(0)}}c_{\la^{(0)}\mu^{(0)}}m_{\mu^{(0)}})\cdots
(\sum_{\mu^{(r-1)}\leq \la^{(r-1)}}c_{\la^{(r-1)}\mu^{(r-1)}}m_{\mu^{(r-1)}}) \\
&=m_{\ula}+\sum_{\umu<\ula}c_{\ula\umu}m_{\umu}, \qquad \mbox{for
some $c_{\ula\umu}\in \mathbb Z$}
\end{align*}
\end{proof}

Let $\ula$ and $\umu$ be colored partitions such that
$\umu\subset\ula$. A \textit{colored skew diagram} $\ula/\umu$ is
the sequence of set-theoretic differences
$\{\lambda^{(i)}-\mu^{(i)}\}_{i\in I}$. A colored maximal connected
component of $\ula/\umu$ is a sequence consisting of maximal
connected components of $\la^{(i)}/\mu^{(i)}$. Each colored
connected component itself is a colored skew diagram. For example
the colored skew Ferrer diagram $\ula /\umu$ =(6,4,2,2)/(4,2) has
three components and is shown by

A \textit{colored horizontal strip} is a sequence of horizontal
strips, i.e., $\{\la^{(i)}/\mu^{(i)}\}_{i\in I}$ (with no two boxes
in the same column) such that $(\la^{(i)})'_k-(\mu^{(i)})'_k\leq 1$.
A \textit{colored vertical strip} is a sequence of vertical strips
$\{\lambda^{(i)} /\mu^{(i)}\}_{i\in I} $, which has no two boxes in
the same row, i.e., ${\la}^{(i)}_k-{\mu}^{(i)}_k\leq 1$.

 A \textit{colored semistandard tableau} $\uT$ of shape $\ula/\umu$ is
 a sequence of semistandard tableau $\{\la^{(i)}/\mu^{(i)}\}_{i\in I}$.
The \textit{content} $\unu$ of a colored tableau $\uT=(T^{(0)},
\ldots, T^{(r-1)})$ is the sequence consisting of the contents
$\nu^{(i)}$ of $\{\la^{(i)}/\mu^{(i)}\}$. For a sequence of
variables $x=(x^{(0)}; \ldots; x^{(r-1)})$ we will write
$x^{\uT}=(x^{(0)})^{T^{(0)}} \cdots (x^{(r-1)})^{T^{(r-1)}}$
associated with the colored tableau $\uT$. Here $x^{(i)}=(x_{i1},
x_{i2}, \ldots)$ for each $i\in I$.

We use
$T_{\ula/\umu, \, \unu}$ to denote the colored skew
tableau of shape $\ula/\umu$ and content $\unu$.
The \textit{Kostka number} $K_{\ula/\umu, \, \unu}$ is the number of
semistandard tableaux of shape $\ula/\umu$ and content $\unu$. For
example, see Table 3. 


\begin{table}[!hpb]
\begin{center}
\begin{tabular}{lll}
\framebox[.28in]{\rule{0in}{.2in}{1}}\framebox[.28in]{\rule{0in}{.2in}{1}}
\framebox[.28in]{\rule{0in}{.2in}{1}}\framebox[.28in]{\rule{0in}{.2in}{4}} &
\framebox[.28in]{\rule{0in}{.2in}{1'}}\framebox[.28in]{\rule{0in}{.2in}{1'}}
\framebox[.28in]{\rule{0in}{.2in}{2'}}\framebox[.28in]{\rule{0in}{.2in}{2'}}&
\framebox[.28in]{\rule{0in}{.2in}{1''}}\framebox[.28in]{\rule{0in}{.2in}{1''}}
\\
\nolinebreak
\framebox[.28in]{\rule{0in}{.2in}{2}}\framebox[.28in]{\rule{0in}{.2in}{3}}
\framebox[.28in]{\rule{0in}{.2in}{3}} &
\framebox[.28in]{\rule{0in}{.2in}{2'}}\framebox[.28in]{\rule{0in}{.2in}{3'}}
\framebox[.28in]{\rule{0in}{.2in}{3'}} &
\framebox[.28in]{\rule{0in}{.2in}{2''}}\\
\framebox[.28in]{\rule{0in}{.2in}{3}}\framebox[.28in]{\rule{0in}{.2in}{4}} &
\framebox[.28in]{\rule{0in}{.2in}{3'}}&
\framebox[.28in]{\rule{0in}{.2in}{3''}}\\
\framebox[.28in]{\rule{0in}{.2in}{5}}\framebox[.28in]{\rule{0in}{.2in}{5}} &
\framebox[.28in]{\rule{0in}{.2in}{5'}}&\\
\framebox[.28in]{\rule{0in}{.2in}{6}} &&
\end{tabular}
\caption{$T_{\ula\umu}$, $\ula=((4,3,2,2,1), (4, 3, 1, 1), (2,1,1))$ and
$\umu=((3,1,3,2,2,1), (2, 3,3,0,1), (2,1,1))$}
\end{center}
\end{table}

\begin{table}[!hpb]
\begin{tabular}{ll}
\framebox[.28in]{\rule{0.2in}{.2in}}\framebox[.28in]{\rule{0.2in}{.2in}}
\framebox[.28in]{\rule{0.2in}{.2in}}\framebox[.28in]{\rule{0.2in}{.2in}}
\framebox[.28in]{\rule{0in}{.2in}{2}}\framebox[.28in]{\rule{0in}{.2in}{2}}&
\framebox[.28in]{\rule{0.2in}{.2in}}\framebox[.28in]{\rule{0.2in}{.2in}}
\framebox[.28in]{\rule{0in}{.2in}{1'}}\framebox[.28in]{\rule{0in}{.2in}{1'}}\\
\framebox[.28in]{\rule{0.2in}{.2in}}\framebox[.28in]{\rule{0.2in}{.2in}}
\framebox[.28in]{\rule{0in}{.2in}{1}}\framebox[.28in]{\rule{0in}{.2in}{1}}
\framebox[.28in]{\rule{0in}{.2in}{5}} &
\framebox[.28in]{\rule{0.2in}{.2in}}
\framebox[.28in]{\rule{0in}{.2in}{1'}}\framebox[.28in]{\rule{0in}{.2in}{2'}}\\
\framebox[.28in]{\rule{0in}{.2in}{2}}\framebox[.28in]{\rule{0in}{.2in}{3}}
\framebox[.28in]{\rule{0in}{.2in}{3}}\framebox[.28in]{\rule{0in}{.2in}{3}} &
\framebox[.28in]{\rule{0in}{.2in}{2'}}\framebox[.28in]{\rule{0in}{.2in}{3'}}\\
\framebox[.28in]{\rule{0in}{.2in}{5}}\framebox[.28in]{\rule{0in}{.2in}{5}}
\framebox[.28in]{\rule{0in}{.2in}{6}}&
\\
\end{tabular}
\caption{tableau $((6,5,4,3), (4,3,2)/((4,2), (2,1))$
with content $((2,3,3,0,3,1), (3, 2, 1))$}
\end{table}
which we will also write as the array

\begin{center}
$\begin{array}{cccccc}
&  &  &  & 2 & 2 \\
&  & 1 & 1 & 5 & \\
2 & 3 & 3 & 3 &  &  \\
5& 5 & 6 &  &  &
\end{array},
\begin{array}{cccccc}
&  &  1' & 1' \\
&  1' & 2' &\\
2' & 3' & &
\end{array}$
\end{center}

The multiplication of two colored Schur functions is
 expressed as a linear combination of colored Schur
functions:
\begin{equation}
s_{\ula}s_{\umu}=\sum_{\unu}c_{\ula\umu}^{\unu}s_{\unu}, \qquad
c_{\ula\umu}^{\unu}\in\mathbb Z_+.
\end{equation}
The structure constants $c_{\ula\umu}^{\unu}$ are called the colored
\textit{Littlewood-Richardson coefficients}. It is known that
$c_{\ula\umu}^{\unu}\neq 0$ unless $|\unu|=|\ula|+|\umu|$ and $\unu\supset
\ula, \umu$.

We define \textit{skew Schur function} $s_{\ula/\umu}\in\LZ$
by
\begin{equation}
s_{\ula/\umu}=\sum_{\unu}c_{\umu\unu}^{\ula}s_{\unu}.
\end{equation}
It is clear that $deg(s_{\ula/\umu})=|\ula|-|\umu|$ and $s_{\ula/\umu}=0$
unless $\umu\subset \ula$, i.e. $s_{\ula/\umu}$ is defined for the skew
Ferrer diagram $\ula/\umu$.

Let $x=(x_1, x_2, \ldots)$ and $y=(y_1, y_2, \ldots)$ be two sets of
variables.
\begin{theorem} \label{schur2}
(a) If the connected components of the colored skew diagram
$\ula/\umu$ are $\uth_j$, then $s_{\ula/\umu}(x)=\prod_j
s_{\uth_j}(x)$. In particular, $s_{\ula/\umu}(x)=0$ if $\umu
\nsubseteq \ula$.

 (b) The skew Schur symmetric function
$s_{\ula}(x, y)$ satisfies
\begin{equation}
s_{\ula/\umu}(x, y)=\sum_{\unu} s_{\ula/\unu}(x)s_{\unu/\umu}(y),
\end{equation}
where the sum runs through all the colored partitions $\unu$ such that
$\ula\supset \unu \supset \umu$.

(c) In general, the skew Schur function $s_{\ula/\umu}(x^{(1)},
\ldots, x^{(n)})$ can be written as
\begin{equation}
s_{\ula/\umu}(x^{(1)}, \ldots, x^{(n)})=\sum_{(\unu)}\prod_{j=1}^n
s_{\unu^{(j)}/\unu^{(j-1)}}(x^{(j)}),
\end{equation}
where the sum runs through all sequences of colored partitions
$(\unu)=(\unu^{(n)}, \ldots, \unu^{(0)})$ such that
$\umu=\unu^{(0)}\subset \unu^{(1)} \subset \ldots \subset\unu^{(n)}
=\ula$.
\end{theorem}

\begin{proof} (a) Suppose the skew diagram $\la^{(i)}/\mu^{(i)}$ has
$m_i$ connected components $\theta^{(i)}_j$, $j=1, \ldots, m_i$.
From Proposition \ref{schur} (a) it follows that
\begin{align*}
s_{\ula/\umu}(x)&=\prod_{i\in I}s_{\la^{(i)}/\mu^{(i)}}(x^{(i)})\\
&=\prod_{i\in I}\prod_{j=1}^{m_i}s_{\theta_{j_i}^{(i)}}(x^{(i)})
=\prod_j^ms_{\uth_j}(x),
\end{align*}
where $m$ is the number of connected components of the colored skew
diagram $\ula/\umu$.

(b) It follows from Proposition \ref{schur} (b) that
\begin{align*}
s_{\ula/\umu}(x, y)&=\prod_{i\in I} s_{\la^{(i)}/\mu^{(i)}}(x^{(i)}, y^{(i)})\\
&=\prod_{i\in I}\left(\sum_{\nu^{(i)}}s_{\la^{(i)}/\nu^{(i)}}(x^{(i)})s_{\nu^{(i)}/\mu^{(i)}}(y^{(i)})\right)\\
&=\sum_{\unu}\prod_{i\in
I}s_{\la^{(i)}/\nu^{(i)}}(x^{(i)})s_{\nu^{(i)}/\mu^{(i)}}(y^{(i)})\\
&=\sum_{\unu} s_{\ula/\unu}(x)s_{\unu/\umu}(y)
\end{align*}
where the sum runs through all colored partitions $\unu$ such that
$\ula\supset \unu \supset \umu$.

Statement (c) follows from (b).
\end{proof}

If there is only one variable in each $x^{(i)}$, and
$x=(x_0; x_1; \cdots; x_{r-1})$, then
the Schur function can be easily
computed.

\begin{theorem} \label{skew2}
(a) We have $s_{\ula}(x)=0$ when $l(\la^{(i)})>1$ for some $i\in I$.
When $\ula/\umu$ is a colored horizontal strip with
$|\la^{(r)}|-|\mu^{(r)}|=a_i$, $i\in I$ and $x=(x_0; x_1; \cdots;
x_{r-1})$, then $s_{\ula/\umu}(x)=x_0^{a_0}\cdots
x_{r-1}^{a_{r-1}}$.

(b) More generally $s_{\ula/\umu}(x_1, \ldots, x_n)=0$ unless
$\ula_i'-\umu_i'\leq n$ for each $i$.
\end{theorem}

\begin{proof} The first identity in Part (a) and Part (b) are a direct
consequence of Proposition \ref{skew}. Now assume
$\ula/\umu$ is a colored horizontal strip with
$|\la^{(r)}|-|\mu^{(r)}|=a_i$, $i\in I$, then
using Proposition \ref{skew} (a) we have

\begin{align*}
s_{\ula/\umu}(x)&=\prod_{i\in I}s_{\la^{(i)}/\mu^{(i)}}(x_i)\\
&=\prod_{i\in I}x_i^{|\la^{(r)}-\mu^{(r)}|}=\prod_{i\in I}x_i^{a_i},
\end{align*}
as required.

\end{proof}

Combining Propositions \ref{skew2} and \ref{schur2} one easily gets
the following result.


\begin{proposition} The Schur function $s_{\ula/\umu}(x)$ can be
expressed as a summation of all monomials $x^{\uT}$ attached to Young
tableaux of shape $\ula/\umu$:
\begin{equation}
s_{\ula/\umu}(x)=\sum_{\uT} x^{\uT},
\end{equation}
where the sum runs through all colored semistandard skew tableaux of
shape $\ula/\umu$.
\end{proposition}
\begin{proof} In Proposition \ref{skew}(c) we take $x^{(i)}=x_i$ and using
Proposition \ref{schur} it follows that
\begin{align*}
s_{\la/\mu}(x_1, \ldots, x_n)&=\sum_{(\nu)}\prod_{i=1}^n
s_{\nu^{(i)}/\nu^{(i-1)}}(x_i)\\
&=\sum_{(\nu)}\prod_{i=1}^n x_i^{|\nu^{(i)}|-|\nu^{(i-1)}|},
\end{align*}
where the sum runs through all sequences of partitions
$(\nu)=(\nu^{(n)}, \ldots, \nu^{(0)})$ such that
$\mu=\nu^{(0)}\subset \nu^{(1)} \subset \ldots \subset\nu^{(n)}
=\la$ and $\nu^{(i)}/\nu^{(i-1)}$ are horizontal strips. For each
horizontal strip $\nu^{(i)}/\nu^{(i-1)}$, we assign $i$ in the
diagram $\la/\mu$ where the horizontal strip occupies and thus
obtain a semistandard Young tableau of shape $\la/\mu$ with content
$T=(|\nu^{(1)}|-|\nu^{(0)}|, |\nu^{(2)}|-|\nu^{(1)}|, \ldots,
|\nu^{(n)}|-|\nu^{(n-1)}|)$. In this way the sequence of such
partitions gives rise a colored semistandard Young tableau of shape
$\ula/\umu$ with content $\uT$. Conversely given a colored semistandard Young
tableau of shape $\ula/\umu$ with content $\uT$ $(|\uT|=|\ula|-|\umu|)$, let
$\nu^{(i)}$ be the subdiagram consisting of the diagram $\mu$
together with entries numbered $1, \ldots, i$. Then the skew diagram
$\nu^{(i)}/\nu^{(i-1)}$ is horizontal and gives rise a sequence of
partitions $(\nu)=(\nu^{(n)}, \ldots, \nu^{(0)})$ such that
$\mu=\nu^{(0)}\subset \nu^{(1)} \subset \ldots \subset\nu^{(n)}
=\la$ and $\nu^{(i)}/\nu^{(i-1)}$ are horizontal strips. Therefore
we have
$$
s_{\ula/\umu}(x)=\sum_{\uT}x^{\uT},
$$
where $\uT$ runs through Young tableaux of shape $\ula/\umu$.
\end{proof}

\begin{corollary} We have $\displaystyle s_{\la/\mu}=\sum_{\nu}K_{\la/\mu,
\nu}m_{\nu}$.
\end{corollary}


\section{Wreath Product of symmetric functions} \label{S:wreath}
 Let $G$ be
a finite group. We denote by $G^*$  its set of
complex irreducible characters and by $G_*=\{c\}$
conjugacy classes and index them as follows.
\begin{equation}
G^*=\{\ga^{(0)}, \ldots, \ga^{(r-1)}\}, \qquad G_*=\{c_0, \ldots,
c_{r-1}\}
\end{equation}
We also let $\ga^{(i)}_j=\ga^{(i)}(c_j)$ and
$\ga^{(i)}_{j^*}=\overline{\ga^{(i)}(c_j)}=\ga^{(i)}(c_j^{-1})$.

Let $\zeta_c$ be the order of the centralizer of the class $c$, then
$\zeta_c=|G|/|c|$. We also write $\zeta_s=\zeta_{c_s}$The orthogonal
relations of the irreducible characters $\ga^{(i)}$ read
\begin{align}\label{orth1}
&\sum_{s=0}^{r-1}{\zeta_s}^{-1}\ga^{(i)}_s\ga^{(j)}_{s^*}=\delta_{ij},\\
&\sum_{i=0}^{r-1}\ga^{(i)}_s\ga^{(i)}_{t^*}=\delta_{st}\zeta_s.
\label{orth2}
\end{align}

Let $p_{n}(\ga)$ ($\ga\in G^*, n\in \mathbb Z_+$) be the independent
indeterminate over $\mathbb C$ associated with the irreducible
character $\ga$ and let
$$
\Lambda_{\mathbb Q}(G)=\mathbb Q[p_{n}(\ga): n\geq 1, \ga\in
G^*]\simeq \Lambda_{\mathbb Q}^{\otimes r}.
$$
For each partition $\lambda=(\la_1, \ldots, \la_l)$, $\la_1\geq
\ldots\geq \la_l\geq 0$, we denote
$p_{\la}(\ga)=p_{\la_1}(\ga)\cdots p_{\la_l}(\ga)$.
 The $\mathbb Q$-algebra $\Lambda_{\mathbb Q}(G)$ is a graded
algebra with degree given by $deg(p_{n}(\ga))=n$. For each colored
partition $\urh=(\rho(\ga))_{\ga\in G^*}$ we define
\begin{equation}\label{powersum1}
p_{\urh}=\prod_{i=0}^{r-1}p_{\rho(\ga^i)}(\ga^i)
\end{equation}

We define the inner product of $\Lambda_{\mathbb Q}(G)$ to be the
tensor product of that of $\Lambda_{\mathbb Q}$. Thus the inner
product is given by
\begin{equation}
<p_{\ula}, p_{\urh}>=\delta_{\ula,\urh}z_{\ula},
\end{equation}
where $z_{\ula}=\prod_{\ga\in G^*}z_{\la(\ga)}$.

We will extend the inner product to the hermitian space
$\Lambda_{\mathbb C}(G)$ naturally (or sesque-linearly) as follows.
For $f=\sum_{\ula}c_{\ula}p_{\ula}$ and
$g=\sum_{\ula}d_{\ula}p_{\ula}$ we define
\begin{equation}\label{innerprod1}
<f, g>=\sum_{\ula}c_{\ula}\overline{d_{\ula}}z_{\ula}.
\end{equation}
In other words, we require that the power sum $p_n(\ga)$ be
invariant under complex conjugation: $\overline{p_{\ula}}=p_{\ula}$.
We remark that this differs from the inner product in \cite{M} where
another power sum basis is invariant.

We now introduce the second power sum symmetric functions indexed by
colored partitions associated with conjugacy classes of $G$. For
each $c\in G_*$ we define
\begin{equation}\label{powersum2}
p_n(c)=\sum_{\ga\in G^*}\ga(c^{-1})p_n(\ga).
\end{equation}
Using the orthogonality of the characters of $G$ we have
\begin{equation}
p_n(\ga)=\sum_{c\in G_*}\zeta_c^{-1}\ga(c)p_n(c)
\end{equation}
For each partition $\lambda=(\la_1, \ldots, \la_l)$, $\la_1\geq
\ldots\geq \la_l\geq 0$, we denote
$p_{\la}(c)=p_{\la_1}(c)\cdots p_{\la_l}(c)$. For a colored partition $\rho$ we define the power sum
\begin{equation}
P_{\urh}=\prod_{i=0}^{r-1} p_{\rho^{(i)}}(c)
\end{equation}
Then $P_{\urh}$ is of degree $n$ if $\urh$ is a colored partition of
$n$. Clearly the power sum symmetric function $P_{\urh}$ forms a
$\mathbb Q$-linear basis for the space $\Lambda^{\otimes r}_{\mathbb
Q}\simeq \Lambda_{\mathbb Q}(G)$. However the second power sum
basis is not invariant under the complex conjugation.

 Then we arrive to the
following equivalent inner product.

\begin{proposition} \label{innerprod2}
  The induced inner product on the wreath product
of the ring $\Lambda$ is given by
\begin{equation}\label{innerprod2}
<P_{\ula}, P_{\urh}>=\delta_{\ula, \urh}Z_{\ula},
\end{equation}
where $Z_{\ula}=\prod_{c_i\in
G_*}z_{\la^{(i)}}\zeta_{c_i}^{l(\la^{(i)})}$ for $\ula$ with
$\la^{(i)}=\la(c_i)$.
\end{proposition}
\begin{proof} Using definition (\ref{powersum2}) and the orthogonal
relations (\ref{orth1}-\ref{orth2}) we compute that
\begin{align*}
<p_m(c),p_n(c')>&=\sum_{\ga, \ga'} \ga(c^{-1})\ga'(c')<p_m(\ga),
p_n(\ga')>\\
&=m\delta_{m,n}\sum_{\ga}\ga(c^{-1})\ga(c')=m\zeta_{c}\delta_{m,
n}\delta_{c, c'}.
\end{align*}
The general case then follows easily by induction.
\end{proof}

For each colored partition $\ula=(\la^{(0)}, \cdots, \la^{(r-1)})$
we define the wreath product Schur symmetric function by
\begin{equation}
S_{\ula}=s_{\la^{(0)}}(P_{\la^{(0)}})\cdots
s_{\la^{(r-1)}}(P_{\la^{(r-1)}}),
\end{equation}
where $s_{\la^{(i)}}(P_{\la^{(i)}})$ is the Schur function spanned
by the power sum $p_{\mu}(\ga_i)$ with $|\mu|=|\la^{(i)}|$.

As our ring $\Lambda_{\mathbb C}(G)$ is isomorphic to the  $r$-fold
tensor product $\Lambda_{\mathbb C}$, we can carry the results from
the previous section to the case of wreath product symmetric
functions.

Likewise we can define wreath product Schur functions for colored
skew partitions. Let $\ula$ and $\umu$ be colored partitions such
that $\umu\subset\ula$. We define \textit{skew wreah product
Schur function}
$S_{\ula/\umu}\in\LZ$ by
\begin{equation}
S_{\ula/\umu}=\sum_{\unu}c_{\umu\unu}^{\ula}S_{\unu},
\end{equation}
where $c_{\umu\unu}^{\ula}$ is the Littewood-Richardson coefficient.

\begin{proposition} The wreath product Schur functions
form an orthonormal basis in the ring $\Lambda_{\mathbb Q}(G)$.
\end{proposition}
\begin{proof} As the wreath product Schur function
is the product of the usual Schur functions, and moreover the
variables of each Schur fucntions are perpendicular to those of the
other Schur functions by (\ref{innerprod2}), the result follows
immediately.
\end{proof}

 Let $x=(x_1, x_2, \ldots)$ and $y=(y_1, y_2, \ldots)$ be two
sets of variables.
\begin{theorem} \label{schur3}
(a) If the connected components of the colored skew diagram
$\ula/\umu$ are $\uth_j$, then $S_{\ula/\umu}(x)=\prod_j
S_{\uth_j}(x)$. In particular, $S_{\ula/\umu}(x)=0$ if $\mu
\nsubseteq \la$.

 (b) The skew wreath product Schur symmetric function
$S_{\ula}(x, y)$ satisfies
\begin{equation}
S_{\ula/\umu}(x, y)=\sum_{\unu} S_{\ula/\unu}(x)S_{\unu/\umu}(y),
\end{equation}
where the sum runs through all the colored partitions $\unu$ such
that $\ula\supset \unu \supset \umu$.

(c) In general, the skew wreath product Schur function
$S_{\ula/\umu}(x^{(1)}, \ldots, x^{(n)})$ can be written as
\begin{equation}
S_{\ula/\umu}(x^{(1)}, \ldots, x^{(n)})=\sum_{(\unu)}\prod_{j=1}^n
S_{\unu^{(j)}/\unu^{(j-1)}}(x^{(j)}),
\end{equation}
where the sum runs through all sequences of colored partitions
$(\unu)=(\unu^{(n)}, \ldots, \unu^{(0)})$ such that
$\umu=\unu^{(0)}\subset \unu^{(1)} \subset \ldots \subset\unu^{(n)}
=\ula$.
\end{theorem}

If there is only one variable in each $x^{(i)}$, and $x=(x_0; x_1;
\cdots; x_{r-1})$, then the wreath product Schur function can be
easily computed.

\begin{theorem} \label{skew2}
(a) We have $S_{\ula}(x)=0$ when $l(\la^{(i)})>1$ for some $i\in I$.
When $\ula/\umu$ is a colored horizontal strip with
$|\la^{(r)}|-|\mu^{(r)}|=a_i$, $i\in I$ and $x=(x_0; x_1; \cdots;
x_{r-1})$, then $s_{\ula/\umu}(x)=x_0^{a_0}\cdots
x_{r-1}^{a_{r-1}}$.

(b) More generally $S_{\ula/\umu}(x_1, \ldots, x_n)=0$ unless
$\ula_i'-\umu_i'\leq n$ for each $i$.
\end{theorem}

\begin{proposition} The wreath product Schur function $S_{\ula/\umu}(x)$ can be
expressed as a summation of all monomials $x^{\uT}$ attached to
Young tableaux of shape $\ula/\umu$:
\begin{equation}
S_{\ula/\umu}(x)=\sum_{\uT} x^{\uT},
\end{equation}
where the sum runs through all colored semistandard skew tableaux of
shape $\ula/\umu$.
\end{proposition}

\begin{corollary} We have $\displaystyle S_{\ula/\umu}=\sum_{\unu}K_{\ula/\umu,
\unu}m_{\unu}$.
\end{corollary}

Let T be a tableau. We derive the \textit{word} of T  or  $w$(T) by
reading the symbols in T from right to left in successive rows,
starting
with the top row. For instance, if T is the tableau\newline\begin{center} $%
\begin{array}{cccccc}
&  &  &  & 5 & 5 \\
&  & 1 & 1 & 6 & 7 \\
2 & 3 & 3 & 3 & 7 & 8 \\
4 & 4 & 6 & 7 & 8 & 9%
\end{array}\newline
$\end{center} $w$(T) is the word 557611873332987644.

If $\uT$ is a colored tableau, then we define the colored word
$\uw(\uT)$ to the composition consists of the words $w(T^{(i)})$. We
also say that the colored word $\uw(\uT)$ is a lattice permutation
if each word $w(T^{(i)})$ is a lattice permutation.

\begin{theorem} \label{LR}. Let $\ \ula ,\umu $, $\unu $ be colored
partitions. Then $c_{\umu \unu }^{\ula }$ is equal to the number of
colored Young tableau $\uT$ of shape $ \ula /\umu $ and content
$\unu $ such that $w(\uT)$ \ is a lattice permutation. The product
$s_{\mu }s_{\nu }$ is an integral linear combination of Schur
functions:
\begin{equation}
S_{\umu }S_{\unu }=\underset{}%
{\underset{\ula }{\sum }c_{\umu \unu }^{\ula }}S_{\ula}
\end{equation}
or equivalently
\begin{equation}
S_{\ula /\umu }=\underset{}{\underset{\unu }{\sum }c_{\umu \unu
}^{\ula }}S_{\unu }
\end{equation}
\end{theorem}
\begin{proof} It follows from the construction of $S_{\ula/\umu}$
that
\begin{equation}
c_{\umu \unu }^{\ula }=c_{\mu^{(0)} \nu^{(0)} }^{\la^{(0)} }\cdots
c_{\mu^{(r-1)} \nu^{(r-1)} }^{\la^{(r-1)} },
\end{equation}
and each $c_{\mu^{(i)} \nu^{(i)} }^{\la^{(i)} }$ is equal to the
number of Young tableau $T$ of shape $ \la^{(i)}/\mu^{(i)} $ and
content $\nu^{(i)} $ such that $w(T)$ \ is a lattice permutation.
Taking product and noticing that they are independent we get the
result in the statement.
\end{proof}

A special case of the Littlewood-Richardson rule is the Pieri rule.

\begin{corollary}
(a) Let $\mathbf{m}=(m_{0},m_{1},\ldots m_{r-1})$ then we have

$S_{\underline{\lambda}}\left(  P\right)  S_{\mathbf{m}}\left(
P\right) =\underset{\underline{\mu}}{\sum}S_{\underline{\mu}}\left(
P\right)  $

summed over all colored partitions \underline{$\mu$} such that $\mu^{(i)}%
/\lambda^{(i)}$ is a horizontal $m_{i}$ strip.

(b) Let $1^{\mathbf{m}}= \left(1^{m_{0}},1^{m_{1}},\ldots1^{
m_{r-1}}\right)$ then we have $S_{\underline{\lambda}}\left(  P
\right)  S_{1^{\mathbf{m}}}\left( P\right)
=\underset{\underline{\mu}}{\sum}S_{\underline{\mu}}\left( P\right)
$

summed over all partitions \underline{$\mu$} such that $\mu^{(i)}%
/\lambda^{(i)}$ is a vertical $m_{i}$ strip.
\end{corollary}

\section{Characters of wreath products}\label{S:char}

For any partition $\lambda$ of $n$, let $\chi^{\lambda}$ be the
irreducible character of $S_n$ corresponding
 to $\lambda$  and $V^{\lambda}$ be an irreducible $\mathbb{C}S_{n}$ module affording $\chi^{\lambda}$.
 Then the degree $\chi^{\lambda}(1)=V^{\lambda}$ is equal to
$n!\slash h(\lambda)$.
 Let $ \left\{ W_{0}, \ldots , W_{r-1}\right)$ be the complete set of representatives of isomorphism classes
 of irreducible $\mathbb C\Gamma$-modules.  We fix the numbering of $W_{0}, \ldots , W_{r-1}$ and put the $d_{i}=dimW_{i}$
. Then the irreducible $\mathbb{C}\left(\Gamma\thicksim
S_{n}\right)$-modules can be constructed as follows. For
$\underline\lambda=\left(\lambda^{0}, \ldots ,
\lambda^{(r-1)}\right)\in (Y^{(r-1)})_{n},$ we put
\begin{center}
$K^{\underline\lambda}=\Gamma^{n}\left(S_{n_{0}}\times\cdots\times
S_{n_{(r-1)}}\right)$ ,\\ \smallskip
$T^{\underline\lambda}=T^{n_{0}}(W_{0})\otimes_{\mathbb{C}}\ldots\otimes_{\mathbb{C}}T^{n_{(r-1)}}(W_{r-1})$,
\\
 and \\
$V^{\underline\lambda}=V^{\lambda^{0}}\otimes_{\mathbb{C}}\ldots\otimes_{\mathbb{C}}V^{\lambda^{(r-1)}}$\\
\end{center}
where $n_{i}=|\lambda^{i}|$ and $T^{i}(V)$ denotes the $i$-fold
tensor product of a vector space $V$. We regard
$T^{\underline\lambda}$ and $V^{\underline\lambda}$ as $\mathbb{C}K^{\underline\lambda}$-modules by defining\\
\begin{center}
$(g_1, \ldots , g_n; \sigma)(t_{1}\otimes\cdots\otimes t_{n})=g_1t_{\sigma^{-1}(1)}\otimes \cdots\otimes g_nt_{\sigma^{-1}(n)}$,\\
$(g_1, \ldots , g_n; \sigma)(v_{0}\otimes\cdots\otimes
v_{(r-1)})=\sigma_{0}v_{0}\otimes\cdots\otimes
\sigma_{(r-1)}v_{(r-1)}$,
\end{center}
where $g=(g_1, \ldots, g_n)\in \Gamma^n$,
$\sigma=\sigma_{1}\cdots\sigma_{(r-1)}\in S_{n_{0}}\times\cdots
\times S_{n_{(r-1)}}$, $t_{0}\otimes\cdots\otimes t_{(r-1)}\in
T^{\underline\lambda}$, and $v_{0}\otimes\cdots\otimes v_{(r-1)}\in
V^{\underline\lambda}$ :
\begin{center}
$W^{\underline\lambda}=\mathbb{C}(\Gamma \thicksim
S_{n})\otimes_{\mathbb{C}K^{\underline\lambda}}(T^{\underline\lambda}\otimes_{\mathbb{C}}V^{\underline\lambda})$.

\end{center}
Then \begin{center}
$dimW^{\underline\lambda}=n!\prod\limits_{i=0}^{r-1}\frac{d_{i}^{|\lambda^{(i)}|}}{h(\lambda^{(i)})}$.

\end{center}
\begin{theorem}\cite{O} For each colored partition $\ula$
the $\mathbb{C}(\Gamma \thicksim S_{n})$-module
$W^{\underline\lambda}$ constructed above is irreducible.  These
$W^{\underline\lambda}$ are pairwise non-isomorphic and exhaust the
isomorphism classes of irreducible $\mathbb{C}(\Gamma \thicksim
S_{n})$ modules.\end{theorem}

Let $f$ be a class function of a finite group $G$, the Frobenius
characteristic $ch: F(G)\longrightarrow \Lambda(G)$, a ring of
symmetric functions associated to $G$,  is defined by
\begin{equation}
ch(f)=\sum_{c}\frac{f(c)}{z_c}p(c)
\end{equation}
where the power sum symmetric function $p(c)$ is associated with the
conjugacy class $c$, and $z_c$ is the order of the conjugacy class
$c$. In all examples of $G$ considered in this paper, the conjugacy
class $c$ is indexed by partitions or colored partitions, so the
meaning of $p(c)$ is self-clear from the context. For example, if
$G=S_n$, then $p(c)$ is the usual power sum symmetric function
corresponding to the cycle-type of the class $c$ and $z_c=z_{\rho}$
with $\rho=\rho(c)$ the cycle type of the class $c$; if $G
=\Gamma\sim S_n$, the conjugacy classes are indexed by colored
partitions $\urh$, and $p(c)=p_{\urh}(c)$, the second type of power
sum symmetric function introduced in Eq. (\ref{powersum2})  and
$z_c=Z_{\urh}$ (cf. Eq. (\ref{innerprod2} )).

The characteristic map enjoys the nice property being an isometric
isomorphism from the ring of class functions to that of symmetric
functions (cf. \cite{M, FJW}).  When $G=S_n$, Frobenius' formula
says that $ch(\chi_{\la})=s_{\la}$, the Schur symmetric function
associated with the partition $\la$. In the case of the wreath
product group $G\sim S_n$ the generalized Frobenius-type formula and
the character table was obtained by Specht \cite{Sp} and recounted
by Macdonald \cite{M}. The idea is to compute the special case of
one-part partition and then extended to general partition using the
isomorphism of $ch$.

\begin{theorem}\cite{Sp, M}
For a family of partitions
$\underline\rho=(\rho^{(0)},\rho^{(1)},\ldots)$ we have
\begin{center}
$S_{\underline\lambda}=\sum\limits_{\underline\rho}Z_{\underline\rho}^{-1}\chi_{\underline\rho}P_{\underline\rho}$
\end{center} or equivalently\begin{center}
$\chi_{\underline\rho}^{\underline\lambda}=\left<
S_{\underline\lambda},P_{\underline\rho}\right>$
\end{center}\end{theorem}
The simplest example of this theory is when $r=2$. We then have the
hyperoctahedral group, which is the semi-direct product
$H_{n}=\mathbb{Z}_{2}^{n}\rtimes S_{n}$, where $\mathbb{Z}_{2}=\{
-1,1\}$ is considered as the multipicative group of two elements.

\bibliographystyle{amsalpha}

\end{document}